\documentclass[11pt,reqno]{amsart}

\usepackage{amscd,amssymb,amsmath,amsthm}
\usepackage[arrow,matrix]{xy}
\usepackage{graphicx}
\usepackage{epstopdf}
\usepackage{color}
\newtheorem{thm}{Theorem}
\newtheorem{defn}{Definition}
\newtheorem{lemma}{Lemma}
\newtheorem{pro}{Proposition}
\newtheorem{rk}{Remark}
\newtheorem{cor}{Corollary}

\numberwithin{equation}{section} 

\def\r{\rho}

\def\Q{\mathbb Q}

\def\N{\mathbb N}
\def\Z{\mathbb{Z}}

\begin{document}

\title[Ergodicity and periodic orbits]{Ergodicity and periodic orbits of $p$-adic $(1,2)$-rational dynamical systems with two fixed points}

\author{I.A. Sattarov, E.T. Aliev}

 \address{I.\ A.\ Sattarov \\ Namangan Satate University,
316, Uychi str., 160100, Namangan, Uzbekistan.} \email
{sattarovi-a@yandex.ru}

 \address{E.\ T.\ Aliev \\ Namangan Institute of Engineering Technology,
7, Kosonsoy str., 160115, Namangan, Uzbekistan.} \email
{aliev-erkinjon@mail.ru}

\begin{abstract}

We consider $(1,2)$-rational functions given on the field of $p$-adic numbers  $\Q_p$. In general, such a function has four parameters. We study the case when such a function has two fixed points and show that when there are two fixed points then $(1,2)$-rational function is conjugate to a two-parametric $(1,2)$-rational function. Depending on these two parameters we determine type of the fixed points, find Siegel disks and the basin of attraction of the fixed points. Moreover, we classify invariant sets and study ergodicity properties of the function on each invariant set. We describe 2- and 3-periodic orbits of the  $p$-adic dynamical systems generated by the two-parametric $(1,2)$-rational functions.

\end{abstract}

 \maketitle

\section{Introduction and preliminaries}

 A function is called a $(n,m)$-rational function if and only if it can be written in the form $f(x)={P_n(x)\over Q_m(x)}$, where
$P_n(x)$ and $Q_m(x)$ are polynomial functions with degree $n$ and
$m$ respectively, $Q_m(x)$ is non zero polynomial.

In this paper we study dynamical systems generated by a (1.2-)rational function. Our investigations based on methods of \cite{ARS}, \cite{AS}, \cite{RS3}-\cite{S}. For motivations of the study see
\cite{AKK}, \cite{A}-\cite{FF}, \cite{MK2}-\cite{UF} and the references therein.

Let us give main definitions.
Let $\Q$ be the field of rational numbers. The greatest common
divisor of the positive integers $n$ and $m$ is denotes by
$(n,m)$. Every rational number $x\neq 0$ can be represented in the
form $x=p^r\frac{n}{m}$, where $r,n\in\mathbb{Z}$, $m$ is a
positive integer, $(p,n)=1$, $(p,m)=1$ and $p$ is a fixed prime
number.

The $p$-adic norm of $x\in \Q$ is given by
$$
|x|_p=\left\{
\begin{array}{ll}
p^{-r}, & \ \textrm{ for $x\neq 0$},\\[2mm]
0, &\ \textrm{ for $x=0$}.\\
\end{array}
\right.
$$
It has the following properties:

1) $|x|_p\geq 0$ and $|x|_p=0$ if and only if $x=0$,

2) $|xy|_p=|x|_p|y|_p$,

3) the strong triangle inequality
$$
|x+y|_p\leq\max\{|x|_p,|y|_p\},
$$

3.1) if $|x|_p\neq |y|_p$ then $|x+y|_p=\max\{|x|_p,|y|_p\}$,

3.2) if $|x|_p=|y|_p$ then $|x+y|_p\leq |x|_p$.

The completion of $\Q$ with  respect to $p$-adic norm defines the
$p$-adic field which is denoted by $\Q_p$ (see \cite{Ko}).

For any $a\in\Q_p$ and
$r>0$ denote
$$
U_r(a)=\{x\in\Q_p : |x-a|_p<r\},\ \ V_r(a)=\{x\in\Q_p :
|x-a|_p\leq r\},
$$
$$
S_r(a)=\{x\in\Q_p : |x-a|_p= r\}.
$$

A function $f:U_r(a)\to\Q_p$ is said to be {\it analytic} if it
can be represented by
$$
f(x)=\sum_{n=0}^{\infty}f_n(x-a)^n, \ \ \ f_n\in \Q_p,
$$ which converges uniformly on the ball $U_r(a)$.

Now let $f:U\to U$ be an analytic function. Denote
$f^n(x)=\underbrace{f\circ\dots\circ f}_n(x)$.

If $f(x_0)=x_0$ then $x_0$
is called a {\it fixed point}. The set of all fixed points of $f$
is denoted by Fix$(f)$. A fixed point $x_0$ is called an {\it
attractor} if there exists a neighborhood $U(x_0)$ of $x_0$ such
that for all points $x\in U(x_0)$ it holds
$\lim\limits_{n\to\infty}f^n(x)=x_0$. If $x_0$ is an attractor
then its {\it basin of attraction} is
$$
A(x_0)=\{x\in \Q_p :\ f^n(x)\to x_0, \ n\to\infty\}.
$$
A fixed point $x_0$ is called {\it repeller} if there  exists a
neighborhood $U(x_0)$ of $x_0$ such that $|f(x)-x_0|_p>|x-x_0|_p$
for $x\in U(x_0)$, $x\neq x_0$.

Let $x_0$ be a fixed point of a
function $f(x)$.
Put $\lambda=f'(x_0)$. The point $x_0$ is attractive if $0<|\lambda|_p < 1$, {\it indifferent} if $|\lambda|_p = 1$,
and repelling if $|\lambda|_p > 1$.

The ball $U_r(x_0)$ is said to
be a {\it Siegel disk} if each sphere $S_{\r}(x_0)$, $\r<r$ is an
invariant sphere of $f(x)$, i.e. if $x\in S_{\r}(x_0)$ then all
iterated points $f^n(x)\in S_{\r}(x_0)$ for all $n=1,2\dots$.  The
union of all Siegel desks with the center at $x_0$ is said to {\it
a maximum Siegel disk} and is denoted by $SI(x_0)$.

Let $f:A\rightarrow A$ and $g:B\rightarrow B$ be two maps. $f$ and $g$ are said
to be {\it topologically conjugate} if there exists a homeomorphism $h: A \rightarrow B$ such
that, $h \circ f = g \circ h$. The homeomorphism $h$ is called a {\it topological conjugacy}.
Mappings which are topologically conjugate are completely equivalent in
terms of their dynamics.

In this paper we consider $(1,2)$-rational function $f:\Q_p\to\Q_p$ defined by
\begin{equation}\label{fa}
f(x)=\frac{ax+b}{x^2+cx+d}, \ \  x\neq \hat x_{1,2}=\frac{-c\pm\sqrt{c^2-4d}}{2}
\end{equation}
where the parameters of the function satisfy the following conditions
$$
 a\neq 0, \ \  a,b,c,d,\sqrt{c^2-4d}\in \Q_p.
$$


We study $p$-adic dynamical systems generated by the rational
function (\ref{fa}).
The equation $f(x)=x$ for fixed points of the function (\ref{fa}) is equivalent to  the equation
\begin{equation}\label{ce}
x^3+cx^2+(d-a)x-b=0.
\end{equation}
The equation (\ref{ce}) may have
three solutions with one of the following relations:

(i) one solution having multiplicity three;

(ii) two solutions, one of which has multiplicity two;

(iii) three distinct solutions.

\begin{rk} Since the behavior of dynamical system depends on the set of fixed points,
each of the above mentioned case (i)-(iii) has its own character of dynamics.  In \cite{RSY} the case (i) was considered.
In this paper we consider the case (ii), i.e., we investigate the behavior of the trajectories of an
arbitrary $(1,2)$-rational dynamical system in $\Q_p$ when  there are two fixed points for $f$. The case (iii) will be considered in a separate
paper.
\end{rk}

The paper is organized as follows. In Section 2 under some assumptions we show that four-parametric function (\ref{fa}) is conjugate to a two-parametric (1,2)-rational function. In Section 3 we study the $p$-adic dynamics generated by the two-parametric  function and give Siegel disks, the basin of attractions and classification of all invariant sets.
In Section 4 we investigate ergodicity of this dynamical systems on invariant sets.
 In Section 5 we describe 2- and 3-periodic orbits.

\section{A function conjugate to (\ref{fa})}

Denote by $x_1$ and $x_2$ the two solutions of the equation (\ref{ce}), where $x_2$ has multiplicity two.
Then we have $x^3+cx^2+(d-a)x-b=(x-x_1)(x-x_2)^2$ and

\begin{equation}\label{f2}\left\{\begin{array}{lll}
x_1+2x_2=-c\\[2mm]
x^2_2+2x_1x_2=d-a\\[2mm]
x_1x^2_2=b.
\end{array} \right.\end{equation}

Let homeomorphism $h:\mathbb Q_p\rightarrow\mathbb Q_p$ be defined by $h(t)=t+x_2$. We note that, the function
$f$ is topologically conjugate to function $h^{-1}\circ f\circ h$. We have

\begin{equation}\label{fh}
(h^{-1}\circ f\circ h)(t)=\frac{-x_2t^2+Bt}{t^2+Dt+B},
\end{equation} where $B=x^2_2+cx_2+d$ and $D=2x_2+c$.

In \cite{RS3} the case $x_2\neq0$ is studied.

Thus in this paper we consider the case $x_2=0$ in (\ref{fh}). If $x_2=0$, then $B=d=a$ and $D=c$. Thus we have the following proposition
\begin{pro}\label{p1}
Any (1,2)-rational function having two distinct fixed points is topologically conjugate to one of the following functions
$$
f(x)=\frac{ax^2+bx}{x^2+cx+b}, \ \ ab(a-c)\neq0, \ \  a,b,c\in \mathbb Q_p,
$$ and
\begin{equation}\label{fc}
f(x)=\frac{ax}{x^2+cx+a}, \ \ ac\neq0, \ \ a,c,\in \mathbb Q_p.
\end{equation}
where  $x\neq \hat x_{1,2}=\frac{-c\pm\sqrt{c^2-4a}}{2}$.

\end{pro}

We study the dynamical system $(\mathbb Q_p, f)$ with $f$ given by (\ref{fc}).

\section{$p$-Adic dynamics of (\ref{fc})}

Note that, the function (\ref{fc}) has two fixed points $x_1=0$ and $x_2=-c$. We have $$f'(x_1)=1 \ \ {\rm and} \ \ f'(x_2)=1-{{c^2}\over{a}}.$$
Thus, the point $x_1$ is an indifferent point for (\ref{fc}), i.e., $x_1$ is a center of some Siegel disk $SI(x_1)$.
In this section we determine the character of the fixed point $x_2$ for each cases.
Then we find Siegel disk or basin of attraction of
the fixed point $x_2$, when $x_2$ is indifferent or attractive, respectively.
In the case where $x_2$ is repelling, we find open ball $U_r(x_2)$, such that the inequality
$|f(x)-x_2|_p>|x-x_2|_p$ holds for all $x\in U_r(x_2)$. Moreover, we study a relation
between the sets $SI(x_1)$ and $SI(x_2)$ when $x_2$ is an indifferent.

For any $x\in \Q_p$, $x\ne \hat x_{1,2}$, by simple calculations we
get
\begin{equation}\label{f2}
    |f(x)|_p=|x|_p\ \cdot{{|a|_p}\over {|x-\hat x_1|_p|x-\hat x_2|_p}}.
\end{equation}

Denote
$$
\mathcal P=\{x\in \Q_p: \exists n\in \N\cup\{0\}, f^n(x)\in\{\hat x_1, \hat x_2\}\},
$$
\begin{equation}\label{0ab}\alpha=\min\{|\hat x_1|_p, |\hat x_2|_p\} \ \ {\rm and} \ \ \beta=\max\{|\hat x_1|_p, |\hat x_2|_p\}.
\end{equation}

Since $\hat x_1+\hat x_2=-c$, we have $|c|_p\leq\alpha$ for $\alpha=\beta$ and $|c|_p=\beta$ for $\alpha<\beta$. Also, since $\hat x_1\hat x_2=a$, we have $|a|_p=\alpha\beta$.

 \begin{thm}\label{t1} The $p$-adic dynamical system generated by the
 function (\ref{fc}) has the following properties:
\begin{itemize}
\item[1.] $SI(x_1)=U_{\alpha}(0)$.
\item[2.] If $|c|_p<\alpha=\beta$, then $x_2$ is indifferent fixed point for (\ref{fc}) and $$SI(x_2)=SI(x_1).$$
\item[3.] If $|c|_p=\alpha=\beta$ and $|a-c^2|_p=\alpha^2$, then $x_2$ is indifferent fixed point for (\ref{fc}) and $$SI(x_2)=U_{\alpha}(x_2), \ \ SI(x_2)\cap SI(x_1)=\emptyset.$$
\item[4.] If $|c|_p=\alpha=\beta$ and $|a-c^2|_p<\alpha^2$, then $x_2$ is attractive fixed point for (\ref{fc}) and $$A(x_2)=U_{\alpha}(x_2)\subset S_{\alpha}(0).$$
\item[5.] If $\alpha<\beta$, then $x_2\in S_{\beta}(0)$ is repelling fixed point for (\ref{fc}) and the inequality
$|f(x)-x_2|_p>|x-x_2|_p$ holds for all $x\in U_{\beta}(x_2)$, $x\neq x_2$.
\end{itemize}
\end{thm}
\begin{proof} 1. Let $x\in S_r(x_1)$, i.e., $|x|_p=r$. Then, from the equalities (3.1), (3.2) and the properties of the $p$-adic norm, we have the following
$$
|f(x)|_p=\left\{\begin{array}{lll}
r,  \ \ \ \ \ \ \mbox{if} \ \ r<\alpha,\\[2mm]
\geq\alpha, \ \ \ \ \ \mbox{if} \ \ \alpha\leq r\leq \beta,\\[2mm]
{{|a|_p}\over r}, \ \ \ \mbox{if} \ \ r>\beta.
\end{array}
\right.
$$
From this equality, $f(S_r(x_1))\subset S_r(x_1)$ for arbitrary $r<\alpha$, i.e. we have  $SI(x_1)=U_{\alpha}(0)$.

2. Note that $|a|_p=\alpha\beta$. If $|c|_p<\alpha=\beta$, then $|f'(x_2)|_p=\left|1-{{c^2}\over a}\right|_p=1.$ From this
$x_2$ is indifferent fixed point.
Let $x\in S_r(x_2)$, i.e., $|x-x_2|_p=r$. Then from the equality
\begin{equation}\label{modx2}
|f(x)-x_2|_p=|x-x_2|_p\cdot{{|x_2(x-x_2)+(x_2^2-a)|_p}\over{|(x-x_2)+\hat x_1|_p|(x-x_2)+\hat x_2|_p}}
 \end{equation}
 we have $|f(x)-x_2|_p=r$ for all $r<\alpha$ and $|f(x)-x_2|_p\geq r$ for $r=\alpha$. Thus,
$f(S_r(x_2))\subset S_r(x_2)$ for arbitrary $r<\alpha$, i.e. we have  $SI(x_2)=U_{\alpha}(x_2)$. In this case, we have $|x_2|_p=|c|_p<\alpha$, so $x_2\in U_{\alpha}(0)=SI(x_1)$. Since these two Siegel disks have the same radii and share a common point, they are the same, i.e., $SI(x_2)=SI(x_1).$

3. If $|c|_p=\alpha=\beta$ and $|a-c^2|_p=\alpha^2$, then $|f'(x_2)|_p=\left|{{a-c^2}\over a}\right|_p=1.$ From this
$x_2$ is indifferent fixed point. As above, from equation (\ref{modx2}) we get $SI(x_2)=U_{\alpha}(x_2)$.
However, in this case $x_2\in S_{\alpha}(0)$, so $SI(x_2)\cap SI(x_1)=\emptyset$.

4. If $|c|_p=\alpha=\beta$ and $|a-c^2|_p<\alpha^2$, then $|f'(x_2)|_p=\left|{{a-c^2}\over a}\right|_p<1.$ From this
$x_2$ is attractive fixed point.  Note that $|x_2|_p=\alpha$. Let $x\in U_{\alpha}(x_2)\subset S_{\alpha}(0)$.
Then from equality (\ref{modx2}) and using the strong triangle inequality of the $p$-adic norm we derive the relation $|f(x)-x_2|_p<|x-x_2|_p$ for all $x\in U_{\alpha}(x_2)$. Similarly, if $x\notin U_{\alpha}(x_2)$, then we have the relation $|f(x)-x_2|_p\geq\alpha$.

Note that, the set of valuations of $p$-adic norm is $\{p^m| \, m\in\Z\}$. Thus, the relation $|f(x)-x_2|_p<|x-x_2|_p$ is equivalent
to the relation $|f(x)-x_2|_p\leq {1\over p}|x-x_2|_p$. This means that the map
$f:U_{\alpha}(x_2)\to U_{\alpha}(x_2)$ is a contraction map. According to the properties of contraction map, we have the equality $A(x_2)=U_{\alpha}(x_2)$.

5. If $\alpha<\beta$, then we have $|x_2|_p=\beta$, i.e., $x_2\in S_{\beta}(0)$. Also, $|f'(x_2)|_p=\left|1-{{c^2}\over a}\right|_p={\beta\over\alpha}>1.$ Let $x\in S_r(x_2)$, i.e., $|x-x_2|_p=r$. Then from the equality (\ref{modx2}) we get
$$
|f(x)-x_2|_p=\left\{\begin{array}{lllll}
{\beta\over\alpha}|x-x_2|_p, \ \ \mbox{if} \ \ r<\alpha,\\[2mm]
\geq\beta, \ \ \ \ \ \ \ \ \ \ \, \mbox{if} \ \ r=\alpha,\\[2mm]
\beta, \ \ \ \ \ \ \ \ \ \ \ \ \ \, \mbox{if} \ \ \alpha<r<\beta,\\[2mm]
\leq\beta, \ \ \ \ \ \ \ \ \ \ \, \mbox{if} \ \ r=\beta,\\[2mm]
\beta, \ \ \ \ \ \ \ \ \ \ \ \ \ \, \mbox{if} \ \ r>\beta.
\end{array}
\right.
$$
From this we conclude that the inequality $|f(x)-x_2|_p>|x-x_2|_p$ is holds for all $x\in U_{\beta}(x_2)$, $x\neq x_2$.
\end{proof}

\begin{cor}\label{cor1}
$\bullet$ The spheres $S_r(x_1)$ is invariant for $f$ if and only if  $r<\alpha$.\\
$\bullet$ The spheres $S_r(x_2)$ is invariant for $f$ if and only if one of the statements holds\\
a) $|c|_p<\alpha=\beta$ and $r<\alpha$;\\
b) $|c|_p=\alpha=\beta$, $|a-c^2|_p=\alpha^2$ and $r<\alpha$.
\end{cor}

\section{Ergodicity of the dynamical systems on invariant spheres}

Recall that an invariant measure is a measure that is preserved by some function.
In ergodic theory of dynamical systems an invariant measure is very important .

Let $G$ be a topological group.
If $G$ is abelian and locally compact, then it is well known \cite{Kan} that it has a nonzero translation-invariant measure $\mu$, which is unique up to scalar. This is called the {\it Haar measure}.

In the field of $p$-adic numbers let $\Sigma$ be the minimal $\sigma$-algebra containing all open and closed (clopen) subsets.

A measure $\mu(V_{\rho})=\rho, \, V_{\rho}\in \Sigma$ is usually called a Haar measure, where $V_{\rho}$ is a ball with radius $\rho$.

However, in some cases, the problem of studying the dynamical system of a function that mapping a compact subset of $\Q_p$ to itself arises. At this time, is needed a measure defined on $\sigma$-algebra with the unit a compact set.
If this compact set has some algebraic structure, then can we look at the natural Haar measure?
If the considered compact set is a ball or a sphere, the answer to this question is positive, which is given as follows in \cite{S1}.

Let $V_r(a)$ be the ball ($S_r(a)$ be the sphere) with the center at the point $a\in\Q_p$ and $\mathcal B$ is the algebra generated by clopen
subsets of $V_r(a)$ ($S_r(a)$). It is known that every element of $\mathcal B$ is a union of
some balls $V_{\rho}(s)\subset V_r(a)$, $s\in V_r(a)$ ($V_{\rho}(s)\subset S_r(a)$, $s\in S_r(a)$).

\begin{thm}\cite{S1}
A measure $\bar\mu:\mathcal B\rightarrow p^{\mathbb Z}$ is a
Haar measure if it is defined by $\bar\mu(V_{\rho}(s))=\rho$ for all $V_{\rho}(s)\in \mathcal B$.
\end{thm}

Also, ergodic theory often deals with ergodic transformations.
Here is the definition:

\begin{defn}\label{erg}\cite{Wal}
Let $T:X\to X$ be a measure-preserving transformation on a measure space $(X, \Sigma, \mu)$, with $\mu(X) = 1$.
Then $T$ is ergodic if for every $E$ in $\Sigma$ with $T^{-1}(E) = E$, either $\mu(E) = 0$ or $\mu(E) = 1$.
\end{defn}

In this section we are interested in ergodicity (with respect to Haar measure)
of the dynamical systems on invariant spheres with the center at the fixed point..

\begin{rk}\label{rk1} Corollary \ref{cor1} in the previous section gives a classification of invariant spheres centered at a fixed point.
Also, in part 2 of Theorem \ref{t1}, it is proved that maximal Siegel discs consisting of union of invariant spheres fall on top of each other. Therefore, the center of invariant spheres is not significant when $|c|_p<\alpha=\beta$. However, when $|c|_p=\alpha=\beta$, it is necessary to consider separately the ergodicity of dynamical systems in invariant spheres with centers $x_1$ and $x_2$.
\end{rk}

For each invariant sphere we consider a measurable space $(S_r(x_i),\mathcal
B)$, here $\mathcal B$ is the algebra generated by closed
subsets of $S_r(x_i)$, $i=1,2$. Every element of $\mathcal B$ is a union of
some balls $V_{\rho}(s)\subset S_r(x_i)$.

A measure $\bar\mu:\mathcal B\rightarrow \mathbb{R}$ is a
\emph{Haar measure} if it is defined by $\bar\mu(V_{\rho}(s))=\rho$ for all $s\in S_r(x_i)$ and $\rho\in p^{\mathbb Z}$ such that $V_{\rho}(s)\subset S_r(x_i)$.

Note that $S_r(x_i)=V_r(x_i)\setminus V_{r\over p}(x_i)$. So, we have
$\bar\mu(S_r(x_i))=r(1-{1\over p})$.

We consider normalized (probability) Haar measure:
$$\mu(V_{\rho}(s))={{\bar\mu(V_{\rho}(s))}\over{\bar\mu(S_r(x_i))}}={{p\rho}\over{(p-1)r}}.$$

\begin{thm}{\label{ab}} Let $S_r(x_i)$, $i=1,2$ be invariant sphere for the function $f$ given by (\ref{fc}).
Then the function $f:S_r(x_i)\to S_r(x_i)$ is an isometry.
\end{thm}
\begin{proof}By the Corollary \ref{cor1}, if the sphere $S_r(x_i)$, $i=1,2$ is invariant for (\ref{fc}), then $r<\alpha$.

Let $i=1$. From relation $x,y\in S_r(x_1)$ we have $|x|_p=|y|_p=r$.
Then, we get the following
\begin{equation}{\label{ab1}}
|f(x)-f(y)|_p=|x-y|_p\cdot\frac{|a|_p|a-xy|_p}{|(x-\hat
x_1)(x-\hat x_2)(y-\hat x_1)(y-\hat x_2)|_p}.
\end{equation}
Note that $|a|_p=\alpha\beta$ and $|x|_p=|y|_p=r<\alpha\leq\beta$. Then,
$$|f(x)-f(y)|_p=|x-y|_p\cdot\frac{\alpha^2\beta^2}{\alpha^2\beta^2}=|x-y|_p.$$
Consequently, the function $f:S_r(x_1)\to S_r(x_1)$ is an isometry.

Let $i=2$. Then by Corollary \ref{cor1} we have two cases. If $|c|_p<\alpha=\beta$ , then by Remark \ref{rk1}, this case overlaps with case $i=1$. If $|c|_p=\alpha=\beta$ and $|a-c^2|_p=\alpha^2$, then by part 3 of Theorem \ref{t1}, we have the relation $S_r(x_2)\subset S_{\alpha}(0)$ for all invariant sphere. So, we have $|x-x_2|_p=r<\alpha$ and $|x|_p=\alpha$ for all $x\in S_{r}(x_2)$.

Let $x,y\in S_{r}(x_2)$. Then
$$
|f(x)-f(y)|_p=|x-y|_p\cdot\frac{|a|_p|(a-x_2^2)+x_2(x_2-y)+y(x_2-x)|_p}{|[(x-x_2)+\hat
x_1][(x-x_2)+\hat x_2][(y-x_2)+\hat x_1][(y-x_2)+\hat x_2]|_p}.
$$
Note that $|a|_p=\alpha^2$, $|x-x_2|_p=|y-x_2|_p=r<\alpha$ and $|a-x_2^2|_p=|a-c^2|_p=\alpha^2$. Then,
$$|f(x)-f(y)|_p=|x-y|_p\cdot\frac{\alpha^4}{\alpha^4}=|x-y|_p.$$
Consequently, the function $f:S_r(x_2)\to S_r(x_2)$ is an isometry.
\end{proof}

\begin{cor}\label{mp}
Let the conditions of the above theorem be satisfied. Then $f:S_r(x_i)\to S_r(x_i)$, $i=1,2$ is a measure-preserving transformation on a measure space $(S_r(x_i), \mathcal B, \mu)$, where $\mu$ is a normalized Haar measure.
\end{cor}

In \cite{S1}, given an important results about the dynamics of isometric maps, and since the function (\ref{fc}) we are considering is also an isometry, the results obtained in \cite{S1} are also relevant for the dynamics of the function (\ref{fc}), i.e., if $S_r(x_i)$, $i=1,2$ is invariant sphere for the function $f$ given by (\ref{fc}), then we have the following:

\begin{itemize}
\item[$\bullet$] The function $f:S_r(x_i)\to S_r(x_i)$, $i=1,2$ is bijection.
\item[$\bullet$] For any initial point $x\in S_r(x_i)$, $i=1,2$ (except fixed point) the orbit
$\{f^n(x)| \, n\in \N\}$ isn't convergent.
\end{itemize}

The result of the following Lemma is given as a condition in \cite{S1}.
Let $S_r(x_i)$, $i=1,2$ be invariant sphere for the function $f$ given by (\ref{fc}), then we denote
$\rho(r,x)=|f(x)-x|_p$ for $x\in S_r(x_i)$.

\begin{lemma} If $r\neq |c|_p$, then for the function $f$ given by (\ref{fc}) the value $\rho(r,x)$ does not depend to $x$.
\end{lemma}
\begin{proof}
We consider all cases in Corollary \ref{cor1}. Let $i=1$. Then $r<\alpha$. By simple calculation we get
$$\rho(r,x)=\left|{{ax}\over{x^2+cx+a}}-x\right|_p=|x|_p^2\cdot{{|x+c|_p}\over{|x-\hat x_1|_p|x-\hat x_2|_p}}=\left\{\begin{array}{ll}
\frac{r^2|c|_p}{\alpha\beta}, \ \ \ \ \  {\rm if} \ \ r<|c|_p,\\[2mm]
\frac{r^3}{\alpha\beta}, \ \ \ \ \ \ \ \ {\rm if} \ \ r>|c|_p.
\end{array}
\right.$$
Let $i=2$. In this case, according to Remark \ref{rk1}, it is sufficient to prove the Lemma when $|c|_p=\alpha=\beta$. So, we have
$r=|x-x_2|_p=|x+c|_p<\alpha$ and

$$\rho(r,x)=\left|{{ax}\over{x^2+cx+a}}-x\right|_p=|x+c|_p\cdot{{|(x+c)-c|_p^2}\over{|(x+c)+\hat x_1|_p|(x+c)+\hat x_2|_p}}=r.$$
\end{proof}

So, we denote $\rho(r)=|f(x)-x|_p$ for all $x\in S_r(x_i)$, $i=1,2$, $r\neq |c|_p$.
In that case, we have the following assertions from \cite{S1}:

\begin{itemize}
\item[$\bullet$] The ball with radius $\rho(r)$ is minimal invariant ball for $f:S_r(x_i)\to S_r(x_i)$, $i=1,2$, $r\neq |c|_p$.

\item[$\bullet$] Let $\mu$ be normalized Haar measure on $S_r(x_i)$. Then
\begin{itemize}
\item[a)] the dynamical system $(S_r(x_i), f, \mu)$ is not ergodic for all $p\geq 3$;
\item[b)] the dynamical system $(S_r(x_i), f, \mu)$ may be ergodic if and only if $r=2\rho(r)$ for $p=2$.
\end{itemize}
\end{itemize}

Let $p=2$. Then according to the above the dynamical system $(S_r(x_2), f, \mu)$ is not ergodic, because $r=\rho(r)$ for $i=2$.

If $i=1$, then $x_1=0$ and we consider the dynamical system $(S_r(0), f, \mu)$.

 Recall  $\mathbb Z_2=\{x\in \Q_2: \, |x|_2\leq 1\}$. So we have $1+2\mathbb Z_2=S_1(0)$.
 The following theorem gives a criterion of ergodicity for the rational functions mapping $S_1(0)$ to itself:

\begin{thm}{\cite{M}}{\label{crit}}
Let $f,g: 1+2\mathbb Z_2\rightarrow 1+2\mathbb Z_2$ be
polynomials whose coefficients are $2$-adic integers.

Set $f(x)=\sum_ia_ix^i$, $g(x)=\sum_ib_ix^i$, and
$$A_1=\sum_{i \, {\rm odd}}a_i, \ \ A_2=\sum_{i \, {\rm even}}a_i, \ \ B_1=\sum_{i \, {\rm odd}}b_i, \ \
B_2=\sum_{i \, {\rm even}}b_i.$$

The rational function $R={f\over g}$ is ergodic if and only if one
of the following situations occurs:

(1) $A_1=1({\rm mod}\,4)$, $A_2=2({\rm mod} \, 4)$, $B_1=0({\rm mod} \, 4)$ and
$B_2=1({\rm mod} \, 4)$.

(2) $A_1=3({\rm mod} 4)$, $A_2=2({\rm mod} \, 4)$, $B_1=0({\rm mod} \, 4)$ and
$B_2=3({\rm mod} \, 4)$.

(3) $A_1=1({\rm mod} \, 4)$, $A_2=0({\rm mod} \, 4)$, $B_1=2({\rm mod} \, 4)$ and
$B_2=1({\rm mod} \, 4)$.

(4) $A_1=3({\rm mod} \, 4)$, $A_2=0({\rm mod} \, 4)$, $B_1=2({\rm mod} \, 4)$ and
$B_2=3({\rm mod} \, 4)$.

(5) One of the previous cases with $f$ and $g$ interchanged.
\end{thm}

Consider $x=g(t)=r^{-1}t$ for $t\in S_1(0)$, then $g^{-1}\circ f\circ g:S_1(0)\rightarrow S_1(0)$.
Let $\mathcal B$ (resp. $\mathcal B_1$) be the algebra generated by closed
subsets of $S_{r}(0)$ (resp. $S_1(0)$), and $\mu$ (resp. $\mu_1$)
be normalized Haar measure on $\mathcal B$ (resp. $\mathcal B_1$).

\begin{thm}{\cite{RS2}}{\label{erg1}} The dynamical system
$(S_{r}(0), \, f, \, \mu)$ is ergodic if and only if\\
$(S_1(0), \, g^{-1}\circ f\circ g, \, \mu_1)$ is ergodic.
\end{thm}

Now using the above mentioned results for
(\ref{fc}) when $p=2$ and
we prove the following

\begin{thm} Let $p=2$. Then the dynamical system $(S_r(0), f, \mu)$ is ergodic if and only if $|c|_2=\beta$ and $r=\frac{\alpha}{2}$.
\end{thm}

\begin{proof} Let $r=2^l$, $\alpha=2^m$, $\beta=2^k$ and $|c|_2=2^q$. Since $\alpha\leq\beta$ we have $m\leq k$. Also, since $c=-\hat x_1-\hat x_2$ and $a=\hat x_1\hat x_2$ we have $q\leq k$ and $|a|_2=2^{m+k}$.

Note that the sphere $S_{2^l}(0)$ is invariant for $f$ iff $l<m$. We consider the function $g: S_1(0)\to S_r(0)$ defined by $x=g(t)=2^{-l}t$. Note that the function\\ $g^{-1}(f(g(t))):S_1(0)\rightarrow S_1(0)$ has the following form
\begin{equation}\label{fg}
g^{-1}(f(g(t)))=\frac{t}{\frac{2^{-2l}}{a}t^2+\frac{2^{-l}c}{a}t+1},
\end{equation}
for the function $f$ given by (\ref{fc}). Note that $k,l,m,q\in\mathbb Z$, $l<m\leq k$ and $q\leq k$. So we have the inequalities
$l-m\leq -1$ and $l-k\leq -1$.
In (\ref{fg}) we can easily see the following $$\left|\frac{2^{-2l}}{a}t^2\right|_2=2^{2l-(m+k)}\leq 2^{-2}, \ \ \ \ \left|\frac{2^{-l}c}{a}t\right|_2=2^{l+q-(m+k)}\leq 2^{-1}.$$

Consequently,
$$t=:\gamma_1(t), \ \ \mbox{is such that} \ \ \gamma_1:1+2\mathbb Z_2\rightarrow 1+2\mathbb Z_2$$
and
$$\frac{2^{-2l}}{a}t^2+\frac{2^{-l}c}{a}t+1=:\gamma_2(t) \ \ \mbox{is such that} \ \ \gamma_2:1+2\mathbb Z_2\rightarrow 1+2\mathbb Z_2.$$

Hence the function (\ref{fg}) satisfies all condition of Theorem \ref{crit}, therefore using this theorem, we get
$$A_1=1, \ \ A_2=0, \ \ B_1=\frac{2^{-l}c}{a} \ \ \mbox{and} \ \ B_2=\frac{2^{-2l}}{a}+1.$$

Moreover,
$$A_1=1({\rm mod}\,4), \ \ A_2=0({\rm mod}\,4), \ \ B_1\in 2^{m+k-(l+q)}(1+2\mathbb Z_2) \ \ \mbox{and} \ \ B_2=1({\rm mod}\,4).$$

By these relations and Theorem \ref{crit} we get $m+k-(l+q)=(m-l)+(k-q)=1.$ Note that $l<m$ and $q\leq k.$ Therefore we conclude that the dynamical system $(S_1(0), \,
g^{-1}\circ f\circ g, \, \mu_1)$ is ergodic if and only if $q=k$ and $l=m-1$, i.e., $|c|_2=\beta$ and $r=\frac{\alpha}{2}.$ Consequently, by Theorem \ref{erg1}, $(S_r(0), f, \mu)$ is ergodic if and only if $|c|_2=\beta$ and $r=\frac{\alpha}{2}.$
\end{proof}

\section{Periodic orbits}

In this section we are interested in periodic trajectories and their characteristics. Since our function is an isometry on an invariant sphere, we get the following result about periodic trajectories from \cite{S1}:

\begin{thm}\label{per}
If the dynamical system $(S_r(x_i), f)$, $i=1,2$ has $n$-periodic orbit $$y_0\to y_1\to ...\to y_n\to y_0,$$
then the following statements hold:
\begin{itemize}
\item[1.] $y_k\in V_{\rho(r)}(y_0)$ for all $k\in\{1, 2, ..., n\}$;
\item[2.] Character of periodic points is indifferent;
\item[3.] If $\rho\leq\rho(r)$, then we have $f(S_{\rho}(y_k))\subset S_{\rho}(y_{k+1})$ for any $k\in\{0,1,...n-1\}$ and\\
$f(S_{\rho}(y_n))\subset S_{\rho}(y_0).$
\end{itemize}
\end{thm}

Now we prove the following theorems about the existence of $2$-periodic and $3$-periodic trajectories:

\begin{thm}\label{per2}
If $\sqrt{c^2-2a}\in\mathbb Q_p$, then the function (\ref{fc}) has unique $2$-periodic orbit $\{t_1, \, t_2\}$,

where $t_{1,2}=-c\pm\sqrt{c^2-2a}$.
\end{thm}

\begin{proof}
We consider the equation
$$\frac{f^2(x)-x}{f(x)-x}=0.$$
Then we obtain the following
$$(x^2+2cx+2a)(x^2+cx+a)=0.$$
Since $x^2+cx+a\neq 0$, we get
$x^2+2cx+2a=0,$ and
$t_{1,2}=-c\pm \sqrt{c^2-2a}.$
\end{proof}

\begin{thm}\label{per3} Let $S_r(x_i)$, $i=1,2$ be invariant sphere for (\ref{fc}) and assume that the parameter $a\in S_r(x_i)$.
Then the function (\ref{fc}) has  $3$-periodic orbit $\left\{a, \, f(a), \, f^2(a)\right\}$ if and only if
\begin{equation}\label{con}
 (a,c)\in \left\{\left(h(q), qh(q)-1\right): q\in \mathbb Q_p\setminus \left\{0, -1, -\frac{2}{3}\right\}, \, |h(q)|_p=r\right\}, \ \
 \mbox{for} \ \ i=1,
\end{equation}
\begin{equation}\label{con1}
(a,c)\in \left\{\left(h(q), qh(q)-1\right): q\in \mathbb Q_p\setminus \left\{0, -1, -\frac{2}{3}\right\}, \, |h(q)(q+1)-1|_p=r\right\}, \ \ \mbox{for} \ \ i=2,
\end{equation}
where $h(q)=\frac{3q^2+2q}{6q^3+11q^2+6q+1}$.
\end{thm}

\begin{proof}
We consider the equation
$$\frac{f^3(x)-x}{f(x)-x}=0.$$
By simplifying this equation, we get the following equation
$$P(x)=x^6+6cx^5+(11c^2+6a)x^4+(6c^3+20ac)x^3+(15ac^2+9a^2)x^2+12a^2cx+3a^3=0.$$
{\it Necessity.}
Let $a\in S_r(x_i)$ be a $3$-periodic point. Then $P(a)=0$ and from this we have the equality
\begin{equation}\label{ac}
a^3+6(c+1)a^2+(11c+9)(c+1)a+3(2c+1)(c+1)^2=0.
\end{equation}
According to equality (\ref{ac}), since $a\neq 0$, we have $c\neq -1$. Denote $$q={{c+1}\over a}.$$
 Then by (\ref{ac}) we get $(6q^3+11q^2+6q+1)a-(3q^2+2q)=0$.

If we denote $$a:=h(q)={{3q^2+2q}\over{6q^3+11q^2+6q+1}},$$ then $c=qh(q)-1$.
Notice that $h(q)$ is undefined at $q=-1$. Applying the conditions that $a(c+1)\neq 0$ we see that $q\neq 0$ and $q\neq -\frac{2}{3}$.

For $i=1$, we have $|a|_p=|h(q)|_p=r$, analogically for $i=2$ we have\\ $|a+c|_p=|h(q)(q+1)-1|_p=r$. Summarizing the above, we get (\ref{con}) and (\ref{con1}).

{\it Sufficiency.} Let conditions (\ref{con}) and (\ref{con1}) be satisfied. Then it is easy to see that $P(a)=0$. Hence, $a\in S_r(x_i)$ is $3$-periodic point for $f$ given by (\ref{fc}).
\end{proof}

\section{Availability of data}
The datasets supporting the conclusions of this article
are included in the article.

\section*{Acknowledgements}
 We thank our supervisor U.A. Rozikov for the useful discussions.

\end{document}